# Q-CALCULUS AND CONVOLUTION TECHNIQUES IN THE STUDY OF Q-RUSCHEWEYEH DERIVATIVES WITH JANOWSKI FUNCTIONS

## K. Marimuthu[1], Nasir Ali[2]

[1]Department of Mathematics, Vel Tech High Tech Dr.Rangarajan Dr.Sukunthala Engineering College, Avadi, Chennai-600062, Tamilnadu, India
[2]Department of Mathematics,COMSATS University Islamabad, Vehari Campus, Pakistan
Email: nasirzawar@gmail.com (Corresponding:Nasir Ali), marimuthu.k@velhightech.com (K. Marimuthu)

ABSTRACT. In this paper, we use concept of q-calculus and technique of convolution to study the q-Ruscheweyeh derivative by the concept of Janowski function, then we define new Subclass of analytic functions. Coefficients Estimates, radii of starlikeness, close to convexity, extreme points and many interesting properties are investigate, obtained and studied.

analytic functions, subordination, q-ruscheweyeh derivative.

## 1. INTRODUCTION

Let $A$ be the class of all analytic functions $f$ defined in the open unit disk $\mathbb{E} = \{z : z \in \mathbb{C}, |z| < 1\}$, normalized by $f(0) = f'(0) - 1 = 0$. Thus for all $f \in A$ are of the form

$$(1.1) \qquad f(z) = z + \sum_{k=2}^{\infty} a_k z^k$$

and let S denote the class of univalent functions in A

Uniformly convex function and starlike function defined and studied by Goodman[1], [2]. Murugusundramoorthy and magesh extended the study of the subclass fixing second coefficient,Extensively studied by Ronning[14], Ma and Minda[13]. The classes of k-uniformly convex and k-uniformly starlike functions introduced by kanas and Wisniowska [11],[12] defined as

$$f \in k - UCV \Leftrightarrow zf'(z) \in k - ST \Leftrightarrow f \in A$$

and

$$Re\left\{\frac{(zf'(z))'}{f'(z)}\right\} > \left|\frac{zf''(z)}{f'(z)}\right|, \ z \in \mathbb{E}, \ k \geq 0.$$

Many authors investigated the properties of $S^*$ and $C$ in multi directions; for a brief study, see[9-19]. let $f, g \in A$. If there exists a Schwartz function $\omega$ analytic in $E$ with $\omega(0) = 0$ and $|\omega(0)| < 1$ for all $z \in \mathbb{E}$ such that $f(z) = g(\omega(z))$, then we say that $f(z)$ is subordinate to $g(z)$ and we write

$$f(z) \prec g(z)$$

Janowski introduced the class $P[A,B]$, using the concept of subordination. The given analytic function with $h(0) = 1$ is belong to the class $P[A,B]$, iff

$$h(z) \prec \frac{1+AZ}{1+BZ}, \quad -1 \leq B < A \leq 1.$$

A function $h(z) \in P[A,B]$. It maps the open unit disk $\mathbb{E}$ onto the domain $\Omega[A,B]$ defined by

(1.2) $$\Omega[A,B] = \omega : \left|\omega - \frac{1-AB}{1-B^2}\right| < \frac{A-B}{1-B^2},$$

this domain represents $D_1 = \frac{A-B}{1-B}$ and $D_2 = \frac{A+B}{1+B}$ with $0 < D_1 < 1 < D_2$.

For $f \in A$ and given by (1.1), $g : g(z) = z + \sum_{k=2}^{\infty} b_k z^k$, Convolution or Hadamard product of f and g defined by

$$(f * g)(z) = \sum_{k=2}^{\infty} a_k b_k z^k, \quad |z| < 1.$$

Recently many researchers involving in the field of geometric function theory because of the use of q- calculus has attracted. In the opening of the last century, q-difference equations studied extensively especially by Mason[6], Jaction[5], Carmichael[4],Trjitzinsky[7], continuous that research many Integral and differential operators can be written in terms of convolution for detail refer to [8]-[18]. To study the quantum groups, fractal and multi-fractal measure, chaotic dynamical systems it helps on q-derivatives and q-integrals, see references [19]-[24]. Ismail et al.[8] generalized the class $S^*$ with concept of q-derivatives and called it $S_q^*$ of q- starlike functions.

For $k$ is an non-negative integer, the $q$-number defined by

(1.3) $$[k,q] = \frac{1-q^k}{1-q}, \quad [0,q] = 0.$$

For $k$ is any non-negative integer, the $q$-number shift factorial is defined by

$$[k,q]! = \begin{cases} 1, & k=0 \\ [1,q][2,q][3,q]\ldots[k,q], & k \in \mathbb{N} \end{cases}$$

Furthermore, the q-generalized Pochhamer symbol for $x > 0$ is given has

$$[k,q]_k = \begin{cases} 1, & k=0 \\ [x,q][x+1,q][x+2,q]\ldots[x+k-1,q], & k \in \mathbb{N} \end{cases}$$

Let the function defined as

(1.4) $$F_{m+1,q}(z) = z + \sum_{k=2}^{\infty} \frac{[m+1,q]_{k-1}}{[k-1,q]!} z^k,$$

where the series convergent absolutely in $\mathbb{E}$. Andrews et al.[19] introduced the q-difference operator related to the q-calculus. Note that when $q \to 1$, $[k-1,q]!$ reduced to the classical definition of factorial. We will assume $q$ to be fixed between $0\,to\,1$, throught this paper.

The q-derivative operator is defined as, for $f \in A$

$$D_q f(z) = \frac{f(qz) - f(z)}{z(q-1)}, z \in \mathbb{E}, q \neq 1.$$

Recently, the Salagean q-differential operator is defined by Govindaraj and Sivasubramanian[20]

$$D_q^0 f(z) = f(z), D_q^1 f(z) = z D_q f(z), \ldots\ldots D_q^k f(z) = z \partial_q (D_q^{k-1} f(z)).$$

The q-Rusheweyeh operator $D_q^m : A \to A$ of order $k \in \mathbb{N}, q \in (0,1)$ then for f given by (1.1) is defined as

(1.5) $$D_q^m f(z) = F_{m+1,q}(z) * f(z) = z + \sum_{k=2}^{\infty} \frac{[m+1,q]_{k-1}}{[k-1,q]!} a_k z^k.$$

Equation (1.5) can be written as

(1.6) $$D_q^m f(z) = \frac{z \partial_q^m (z^{m-1} f(z))}{[m,q]!}, m \in \mathbb{N}_0.$$

Since $\lim_{q \to 1^-} F_{m+1,q}(z) = \frac{z}{(1-z)^{m+1}}$, it follows that $\lim_{q \to 1} D_q^m f(z) = \frac{z}{(1-z)^{m+1}} * f(z) = D^m f(z)$.

We introduced a new subclass $TM_{m,l}(q, \alpha, A, B)$ from the above-cited work.

**Definition 1.1.** A function $f \in A$ is in the subclass $M_{m,l}(q, \alpha, A, B)$ of Subordination condition,

$$\frac{D_q^m f(z)}{D_q^l f(z)} - \alpha \left| \frac{D_q^m f(z)}{D_q^l f(z)} - 1 \right| \prec \frac{1+Az}{1+Bz}$$

where $-1 \leq B < A \leq 1, \alpha \leq 0, m \in \mathbb{N}, l \in \mathbb{N}_0, m > l, q \in (0,1), z \in \mathbb{E}$.

**Definition 1.2.** Let $M$ denote the subclass of functions of $A$ of the form

(1.7) $$f(z) = z - \sum_{k=2}^{\infty} a_k z^k, \ a_k \leq 0.$$

Further, we define the class $TM_{m,l}(q,\alpha,A,B) = M_{m,l}(q,\alpha,A,B) \cap M$.

## 2. Coefficients inequalities

In this section we will prove our main results.

**Theorem 2.1.** *A function of the form* (1.1) *is in the class* $M_{m,l}(q,\alpha,A,B)$ *if:*

$$\sum_{k=2}^{\infty}\{(1+\alpha(1+|B|))([m+1,q]_{k-1}-[l+1,q]_{k-1})+|B[m+1,q]_{k-1}-A[l+1,q]_{k-1}|\}a_k \leq A-B \quad (2.1)$$

*Proof.* It is sufficient to show that: $\left|\frac{p(z)-1}{A-Bp(z)}\right| < 1$ where $p(z) = \frac{D_q^m f(z)}{D_q^l f(z)} - \alpha\left|\frac{D_q^m f(z)}{D_q^l f(z)} - 1\right|$.

We have,

$$\left|\frac{p(z)-1}{A-Bp(z)}\right|$$

$$= \left|\frac{D_q^m f(z) - D_q^l f(z) - \alpha e^{i\theta}|D_q^m f(z) - D_q^l f(z)|}{AD_q^l f(z) - B[D_q^m f(z) - \alpha e^{i\theta}|D_q^m f(z) - D_q^l f(z)|]}\right|$$

$$= \left|\frac{\sum_{k=2}^{\infty}([m+1,q]_{k-1} - [l+1,q]_{k-1})a_k z^k - \alpha e^{i\theta}|\sum_{k=2}^{\infty}([m+1,q]_{k-1} - [l+1,q]_{k-1})a_k z^k|}{(A-B)z - \left\{\sum_{k=2}^{\infty}(B[m+1,q]_{k-1} - A[l+1,q]_{k-1})a_k z^k - B\alpha e^{i\theta}|\sum_{k=2}^{\infty}([m+1,q]_{k-1} - [l+1,q]_{k-1})a_k\right.}\right|$$

$$\leq \frac{\sum_{k=2}^{\infty}([m+1,q]_{k-1} - [l+1,q]_{k-1})|a_k||z|^k - \alpha \sum_{k=2}^{\infty}([m+1,q]_{k-1} - [l+1,q]_{k-1})|a_k||z|^k}{(A-B)|z| - \left\{\sum_{k=2}^{\infty}|B[m+1,q]_{k-1} - A[l+1,q]_{k-1}||a_k||z|^k - B\alpha \sum_{k=2}^{\infty}([m+1,q]_{k-1} - [l+1,q]_{k-1})|a\right.}$$

$$\leq \frac{\sum_{k=2}^{\infty}([m+1,q]_{k-1} - [l+1,q]_{k-1})(1+\alpha)|a_k|}{(A-B) - \sum_{k=2}^{\infty}|B[m+1,q]_{k-1} - A[l+1,q]_{k-1}||a_k| - \alpha|B|\sum_{k=2}^{\infty}([m+1,q]_{k-1} - [l+1,q]_{k-1})|a_k|}$$

This last expression is bounded above by one if:

$$\sum_{k=2}^{\infty}\{(1+\alpha(1+|B|))([m+1,q]_{k-1}-[l+1,q]_{k-1})+|B[m+1,q]_{k-1}-A[l+1,q]_{k-1}|\}a_k \leq A-B$$

Hence, the proof is completed.

□

Theorem(3.2) shows that the condition(3.1) is also necessary for functions $f$ of the form (2.1) to be in the class $TM_{i,j}(q,\alpha,A,B)$.

**Theorem 2.2.** *Let* $f \in T$. *Then,* $TM_{i,j}(q,\alpha,A,B)$ *if and only if*

$$\sum_{k=2}^{\infty}\{(1+\alpha(1+|B|))([m+1,q]_{k-1}-[l+1,q]_{k-1})+|B[m+1,q]_{k-1}-A[l+1,q]_{k-1}|\}a_k \leq A-B$$

*Proof.* Since $TM_{i,j}(q,\alpha,A,B) \subset M_{i,j}(q,\alpha,A,B)$, for functions $f \in TM_{i,j}(q,\alpha,A,B)$, we can write

$$\left|\frac{p(z)-1}{A-Bp(z)}\right| < 1,$$

where

$$p(z) = \frac{D_q^m f(z)}{D_q^l f(z)} - \alpha\left|\frac{D_q^m f(z)}{D_q^l f(z)} - 1\right|.$$

Then

$$\left|\frac{\sum_{k=2}^{\infty}([m+1,q]_{k-1}-[l+1,q]_{k-1})a_k z^k + \alpha e^{i\theta}|\sum_{k=2}^{\infty}([m+1,q]_{k-1}-[l+1,q]_{k-1})a_k z^k|}{(A-B)z + \left\{\sum_{k=2}^{\infty}(B[m+1,q]_{k-1}-A[l+1,q]_{k-1})a_k z^k + B\alpha e^{i\theta}|\sum_{k=2}^{\infty}([m+1,q]_{k-1}-[l+1,q]_{k-1})a_k z^k|\right\}}\right|$$

Since $Re(z) \leq |z|$, then we obtain

$$Re\left(\frac{\sum_{k=2}^{\infty}([m+1,q]_{k-1}-[l+1,q]_{k-1})a_k z^k + \alpha e^{i\theta}|\sum_{k=2}^{\infty}([m+1,q]_{k-1}-[l+1,q]_{k-1})a_k z^k|}{(A-B)z + \left\{\sum_{k=2}^{\infty}(B[m+1,q]_{k-1}-A[l+1,q]_{k-1})a_k z^k + B\alpha e^{i\theta}|\sum_{k=2}^{\infty}([m+1,q]_{k-1}-[l+1,q]_{k-1})a_k z^k\right.}\right.$$

Now choosing z to be real and letting $z \to 1^-$, we obtain

$$\sum_{k=2}^{\infty}\{(1+\alpha(1+|B|))([m+1,q]_{k-1}-[l+1,q]_{k-1})+|B[m+1,q]_{k-1}-A[l+1,q]_{k-1}|\}a_k \leq A-B;$$

or equivalently

$$\sum_{k=2}^{\infty}\{(1+\alpha(1+|B|))([m+1,q]_{k-1}-[l+1,q]_{k-1})+|B[m+1,q]_{k-1}-A[l+1,q]_{k-1}|\}a_k \leq A-B.$$

Hence completes the proof. □

**Corollary 2.1.** *Let $f$ defined by (2.1) be in the class $TM_{i,j}(q,\alpha,A,B)$. Then*

(2.2)
$$a_k \leq \frac{A-B}{\{(1+\alpha(1+|B|))([m+1,q]_{k-1}-[l+1,q]_{k-1})+|B[m+1,q]_{k-1}-A[l+1,q]_{k-1}|\}},\ m,l \geq 1$$

*The result is sharp for the function*

(2.3)
$$f(z) = z - \frac{A-B}{\{(1+\alpha(1+|B|))([m+1,q]_{k-1}-[l+1,q]_{k-1})+|B[m+1,q]_{k-1}-A[l+1,q]_{k-1}|\}}z^2,\ m,l \geq$$

*That is equality can be attained for the function defined in (3.3).*

## 3. Distortion theorems

**Theorem 3.1.** *Let the function $f$ be defined by (2.1) be in the class $TM_{i,j}(q, \alpha, A, B)$. Then*

$$|f(z)| \geq |z| - \frac{A-B}{\{(1+\alpha(1+|B|))([2,q]_{k-1} - [2,q]_{k-1}) + |B[2,q]_{k-1} - A[2,q]_{k-1}|\}}|z|^2,$$

*and*

$$|f(z)| \leq |z| + \frac{A-B}{\{(1+\alpha(1+|B|))([2,q]_{k-1} - [2,q]_{k-1}) + |B[2,q]_{k-1} - A[2,q]_{k-1}|\}}|z|^2$$

*The result is sharp.*

*Proof.* In view of theorem 3.2, consider the function

$$\mu(k) = \{(1+\alpha(1+|B|))([m+1,q]_{k-1} - [l+1,q]_{k-1}) + |B[m+1,q]_{k-1} - A[l+1,q]_{k-1}|\}.$$

then it is clear that it is an increasing function of $m$ and $l$ $(m, l \geq 1)$, therefore

$$\mu(2) \sum_{k=2}^{\infty} |a_k| \leq \sum_{k=2}^{\infty} \mu(k)|a_k| \leq A - B.$$

That is

$$\sum_{k=2}^{\infty} |a_k| \leq \frac{A-B}{\mu(2)}.$$

Thus we have

$$|f(z)| \leq |z| + |z|^2 \sum_{k=2}^{\infty} |a_k|,$$

$$|f(z)| \leq |z| + \frac{A-B}{\{(1+\alpha(1+|B|))([2,q]_{k-1} - [2,q]_{k-1}) + |B[2,q]_{k-1} - A[2,q]_{k-1}|\}}|z|^2.$$

Similarly, we get

$$|f(z)| \geq |z| - \sum_{k=2}^{\infty} |a_k||z|^k \geq |z|^2 \sum_{k=2}^{\infty} |a_k|$$

$$\geq |z| - \frac{A-B}{\{(1+\alpha(1+|B|))([2,q]_{k-1} - [2,q]_{k-1}) + |B[2,q]_{k-1} - A[2,q]_{k-1}|\}}|z|^2.$$

Finally, the equality can be attained for the function

$$(3.1) \quad f(z) = z - \frac{A-B}{\{(1+\alpha(1+|B|))([2,q]_{k-1} - [2,q]_{k-1}) + |B[2,q]_{k-1} - A[2,q]_{k-1}|\}}z^2.$$

At $|z| = r$ and $z = re^{i(2t+1)\pi}$ $(t \in Z)$. This completes the result. □

**Theorem 3.2.** *Let the function $f$ be defined by (2.1) in the class $TM_{i,j}(q, \alpha, A, B)$. Then*

$$|f'(z)| \geq 1 - \frac{A-B}{\{(1+\alpha(1+|B|))([2,q]_{k-1} - [2,q]_{k-1}) + |B[2,q]_{k-1} - A[2,q]_{k-1}|\}}|z|,$$

*and*

$$|f'(z)| \leq 1 + \frac{A-B}{\{(1+\alpha(1+|B|))([2,q]_{k-1} - [2,q]_{k-1}) + |B[2,q]_{k-1} - A[2,q]_{k-1}|\}}|z|.$$

*The result is sharp.*

*Proof.* As $\frac{\mu(k)}{k}$ is an increasing function for $k(k \geq 1)$, in view of theorem (3.2) we have

$$\frac{\mu(2)}{2}\sum_{k=2}^{\infty} k|a_k| \leq \sum_{k=2}^{\infty}\frac{\mu(k)}{k}n|a_k| = \sum_{k=2}^{\infty}\mu(k)|a_k| \leq (A-B),$$

that is

$$\sum_{k=2}^{\infty} k|a_k| \leq \frac{2(A-B)}{\mu(2)}.$$

Thus we have

$$|f'(z)| \leq 1 + |z|\sum_{k=2}^{\infty} k|a_k|,$$

$$|f(z)| \leq 1 + \frac{2(A-B)}{\{(1+\alpha(1+|B|))([2,q]_{k-1} - [2,q]_{k-1}) + |B[2,q]_{k-1} - A[2,q]_{k-1}|\}}|z|.$$

Similarly we get

$$|f(z)| \geq 1 - |z|\sum_{k=2}^{\infty} k|a_k|$$

$$\geq 1 - \frac{2(A-B)}{\{(1+\alpha(1+|B|))([2,q]_{k-1} - [2,q]_{k-1}) + |B[2,q]_{k-1} - A[2,q]_{k-1}|\}}|z|.$$

Finally, we can see that the assertions of the theorem are sharp for the function f(z) defined by (3.4). This completes the proof. □

## 4. Radii of starlikeness, convexity and close-to-convexity

**Theorem 4.1.** *Let the function f defined by (2.1) be in the class $TM_{i,j}(q,\alpha,A,B)$. Then*

*(i) f is starlike of order $\psi(0 \leq \psi < 1)$ in $|z| < r_1$, where*

(4.1)
$$r_1 = \inf_{p\geq 2}\left\{\frac{((1+\alpha(1+|B|))([m+1,q]_{k-1} - [l+1,q]_{k-1}) + |B[m+1,q]_{k-1} - A[l+1,q]_{k-1}|)}{(A-B)}X\left(\frac{1-}{(p-}\right.\right.$$

*(ii) f is convex of order $\psi(0 \leq \psi < 1)$ in $|z| < r_2$, where*

(4.2)
$$r_2 = \inf_{p\geq 2}\left\{\frac{((1+\alpha(1+|B|))([m+1,q]_{k-1} - [l+1,q]_{k-1}) + |B[m+1,q]_{k-1} - A[l+1,q]_{k-1}|)}{(A-B)}X\left(\frac{1-}{p(p}\right.\right.$$

*(iii) f is close convex of order $\psi(0 \leq \psi < 1)$ in $|z| < r_3$, where*

(4.3)
$$r_3 = \inf_{p\geq 2}\left\{\frac{((1+\alpha(1+|B|))([m+1,q]_{k-1} - [l+1,q]_{k-1}) + |B[m+1,q]_{k-1} - A[l+1,q]_{k-1}|)}{(A-B)}X\left(\frac{1-}{p}\right.\right.$$

*Each of these results is sharp for the function f given by (3.3)*

*Proof.* It is sufficient to show that

$$\left|\frac{zf'(z)}{f(z)} - 1\right| \leq 1 - \mu \text{ for } |z| < r_1,$$

where $r_1$ is given by (3.5). Indeed we find from (2.1) that

$$\left|\frac{zf'(z)}{f(z)} - 1\right| \leq \frac{\sum_{k=2}^{\infty}(p-1)a_k|z|^{k-1}}{1 - \sum_{k}^{\infty} a_k|z|^{k-1}}.$$

Thus we have

$$\left|\frac{zf'(z)}{f(z)} - 1\right| \leq 1 - \psi,$$

if and only if

(4.4)
$$\frac{\sum_{k=2}^{\infty}(p-\psi)a_k|z|^{k-1}}{(1-\psi)} \leq 1.$$

From theorem 3.2, the relation (3.8) is true if

$$\left(\frac{p-\psi}{1-\psi}\right)|z|^{k-1} \leq \left\{\frac{((1+\alpha(1+|B|))([m+1,q]_{k-1} - [l+1,q]_{k-1}) + |B[m+1,q]_{k-1} - A[l+1,q]_{k-1}|)}{(A-B)}\right\}$$

That is, if

$$|z| \leq \left\{\frac{((1+\alpha(1+|B|))([m+1,q]_{k-1} - [l+1,q]_{k-1}) + |B[m+1,q]_{k-1} - A[l+1,q]_{k-1}|)}{(A-B)}X\left(\frac{1-\psi}{(p-\psi)}\right)\right\}$$

This implies that

$$r_1 = \inf_{p\geq 2}\left(\frac{\{(1+\alpha(1+|B|))([m+1,q]_{k-1} - [l+1,q]_{k-1}) + |B[m+1,q]_{k-1} - A[l+1,q]_{k-1}|\}}{(A-B)}X\left(\frac{1-}{(p-}\right.\right.$$

This completes the proof of (3.5) To prove (3.6) and (3.7), it is sufficient to show that

$$\left|1 + \frac{zf''(z)}{f'(z)} - 1\right| \leq 1 - \psi \quad (|z| < r_2, 0 \leq \psi < 1),$$

and

$$\left|f'(z) - 1\right| \leq 1 - \psi \quad (|z| < r_3, 0 \leq \psi < 1),$$

respectively.

## 5. Extreme points

**Theorem 5.1.** *Let $f_1(z) = z$, and*

$$f_k(z) = z - \frac{(A-B)}{\{(1+\alpha(1+|B|))([m+1,q]_{k-1} - [l+1,q]_{k-1}) + |B[m+1,q]_{k-1} - A[l+1,q]_{k-1}|\}}z^k, \; k = 2,$$

*Then, $f \in TM_{i,j}(q,\alpha,A,B)$. iff it can be expressed in the following form*

$$f(z) = \sum_{k=1}^{\infty}\eta_k f_k(z),$$

where

$$\eta_k \geq 0, \; \sum_{k=1}^{\infty}\eta_k = 1.$$

*Proof.* Suppose that
$$f(z) = \sum_{k=1}^{\infty} \eta_k f_k(z)$$

$$z - \sum_{k=2}^{\infty} \eta_k \frac{(A-B)}{\{(1+\alpha(1+|B|))([m+1,q]_{k-1} - [l+1,q]_{k-1}) + |B[m+1,q]_{k-1} - A[l+1,q]_{k-1}|\}} z^k.$$

Then, from Theorem 3.2, we have

$$\sum_{k=2}^{\infty} \left[ \frac{\{(1+\alpha(1+|B|))([m+1,q]_{k-1} - [l+1,q]_{k-1}) + |B[m+1,q]_{k-1} - A[l+1,q]_{k-1}|\}(A-B)}{\{(1+\alpha(1+|B|))([m+1,q]_{k-1} - [l+1,q]_{k-1}) + |B[m+1,q]_{k-1} - A[l+1,q]_{k-1}|\}} \eta_p \right]$$

$$= (A-B) \sum_{k=2}^{\infty} \eta_k = (A-B)(1-\eta_1) \leq (A-B).$$

Thus, in view of Theorem 3.2, we find that $TM_{i,j}(q,\alpha,A,B)$. Conversely, let us suppose that $f \in TM_{i,j}(q,\alpha,A,B)$, then

$$a_k \leq \frac{(A-B)}{\{(1+\alpha(1+|B|))([m+1,q]_{k-1} - [l+1,q]_{k-1}) + |B[m+1,q]_{k-1} - A[l+1,q]_{k-1}|\}}$$

By setting

$$\psi_k = \frac{\{(1+\alpha(1+|B|))([m+1,q]_{k-1} - [l+1,q]_{k-1}) + |B[m+1,q]_{k-1} - A[l+1,q]_{k-1}|\}}{(A-B)} a_k$$

with

$$\psi_1 = 1 - \sum_{k=2}^{\infty},$$

$$f(z) = \sum_{k=1}^{\infty} \psi_k f_k(z).$$

This completes the proof.

**Corollary 5.1.** *The extreme points of the class $TM_{i,j}(q,\alpha,A,B)$ are given by*

$$f_1(z) = z,$$

*and*

$$f_k(z) = z - \frac{(A-B)}{\{(1+\alpha(1+|B|))([m+1,q]_{k-1} - [l+1,q]_{k-1}) + |B[m+1,q]_{k-1} - A[l+1,q]_{k-1}|\}} z^k, \quad k = 2,$$

## 6. Integral mean inequalities

**Lemma 6.1.** *If the functions $f$ and $g$ are analyticin $E$ with*

$$f_1(z) \prec g(z),$$

*then for $s > 0$ and $z = re^{i\theta}$, $(0 < r < 1)$,*

(6.1) $$\int_0^{2\pi} |f(z)|^s d\theta \leq \int_0^{2\pi} |g(z)|^s d\theta.$$

*We now make use of lemma to prove the following result.*

**Theorem 6.1.** *Suppose that $TM_{i,j}(q, \alpha, A, B)$, $s > 0, -1 \leq B < A \leq 1, \alpha > 0, m \in N, l \in N_0, m > l$, and $f_2(z)$ is defined by*

$$f_2(z) = z - \frac{(A-B)}{\{(1+\alpha(1+|B|))([2,q]_{k-1} - [2,q]_{k-1}) + |B[2,q]_{k-1} - A[2,q]_{k-1}|\}} z^2,$$

*then for $z = re^{i\theta}$, $(0 < r < 1)$, we have*

$$\int_0^{2\pi} |f(z)|^k d\theta \leq \int_0^{2\pi} |f_2(z)|^k d\theta.$$

*Proof.* For

$$f(z) = z - \sum_{k=2}^\infty a_k z^k, \quad a_k \geq 0,$$

the relation (3.9) is equivalent to proving that

$$\int_0^{2\pi} \left|1 - \sum_{k=2}^\infty a_k z^{k-1}\right|^k d\theta$$

$$\leq \int_0^{2\pi} \left|1 - \frac{(A-B)}{\{(1+\alpha(1+|B|))([2,q]_{k-1} - [2,q]_{k-1}) + |B[2,q]_{k-1} - A[2,q]_{k-1}|\}} z\right|^k d\theta.$$

By applying lemma it suffices to show that

$$1 - \sum_{k=2}^\infty a_k z^{k-1} \prec 1 - \frac{(A-B)}{\{(1+\alpha(1+|B|))([2,q]_{k-1} - [2,q]_{k-1}) + |B[2,q]_{k-1} - A[2,q]_{k-1}|\}} z.$$

By setting

$$1 - \sum_{k=2}^\infty a_k z^{k-1} = 1 - \frac{(A-B)}{\{(1+\alpha(1+|B|))([2,q]_{k-1} - [2,q]_{k-1}) + |B[2,q]_{k-1} - A[2,q]_{k-1}|\}} \omega(z),$$

and using (3.1), we obtain

$$\begin{aligned}|\omega(z)| &= \left|\sum_{k=2}^{\infty} \frac{\{(1+\alpha(1+|B|))([2,q]_{k-1} - [2,q]_{k-1}) + |B[2,q]_{k-1} - A[2,q]_{k-1}|\}}{(A-B)} a_k z^{k-1}\right| \\ &\leq |z| \sum_{k=2}^{\infty} \frac{\{(1+\alpha(1+|B|))([2,q]_{k-1} - [2,q]_{k-1}) + |B[2,q]_{k-1} - A[2,q]_{k-1}|\}}{(A-B)} a_k \\ &\leq \sum_{k=2}^{\infty} \frac{\{(1+\alpha(1+|B|))([m+1,q]_{k-1} - [l+1,q]_{k-1}) + |B[m+1,q]_{k-1} - A[l+1,q]_{k-1}|\}}{(A-B)} a_k \\ &\leq |z| < 1\end{aligned}$$

This completes the proof.